\documentclass{amsart}
\usepackage{url}

\usepackage{graphicx}
\usepackage{datetime}
\usepackage{color}
\usepackage{amsmath,amsfonts,amsthm,amssymb}
\usepackage{epstopdf}
\usepackage[all,cmtip]{xy}
\usepackage{accents}
\newtheorem{theorem}{Theorem}
\newtheorem{proposition}[theorem]{Proposition}

\newtheorem{corollary}[theorem]{Corollary}
\newtheorem{lemma}[theorem]{Lemma}
\newtheorem{remark}[theorem]{Remark}
\numberwithin{equation}{section} \numberwithin{theorem}{section}

\newcommand{\CM}{\mathcal{M}}
\newcommand{\R}{\mathbb{R}}
\newcommand{\e}{\varepsilon}
\newcommand{\N}{\mathbb{N}}

\newcommand{\mc}{\mathcal}

\def\XXint#1#2#3{{\setbox0=\hbox{$#1{#2#3}{\int}$}
     \vcenter{\hbox{$#2#3$}}\kern-.5\wd0}}

\newcommand{\ov}{\overline}

\DeclareMathOperator{\diam}{diam}

\begin{document}

\hyphenation{non-ne-ga-ti-ve}
\title[Tori Gap Theorem and Willmore Flow]{A Gap Theorem for Willmore Tori and an application to the Willmore Flow}  

\author{Andrea Mondino}
\address[A. Mondino]{ ETH, R\"amistrasse 101, Zurich, Switzerland}
\email{andrea.mondino@math.ethz.ch}

\author{Huy The Nguyen}
\address[H.T.Nguyen]{School of Mathematics and Physics\\
The University of Queensland\\
St Lucia, Brisbane, Australia, 4072}
\email{huy.nguyen@maths.uq.edu.au}

\thanks{The first author is supported by the ETH Fellowship; the second author was supported by the Leverhulme Trust.}

\begin{abstract} 
In 1965 Willmore conjectured that the integral of the square of the mean curvature of a torus immersed in $ \R^3$ is at least $2\pi^2$ and attains this minimal value if and only if the torus is a M\"obius transform of the Clifford torus. This was recently proved by Marques and Neves in  \cite{MN}. In this paper, we show for tori there is a gap to the next critical point of the Willmore energy and we discuss an application to  the Willmore flow. We also prove an energy gap from the Clifford torus to surfaces of higher genus.   
\end{abstract}
\maketitle

\section{Introduction}

Let $\Sigma$ be a compact Riemann surface and $ f: \Sigma\rightarrow \R^3$ be a smooth immersion. Then the Willmore energy is defined to be 
\begin{align*}
\mc W (f)= \int |H|^2 d \mu _g 
\end{align*} 
where $g$ is the induced metric, $ d \mu_g$ is the induced area form and $H$ is the mean curvature  (we adopt the convention that $H$ is half of the trace of the second fundamental  form). It is well known that the Willmore energy is invariant under conformal transformations of $\R^3$, the so called M\"obius transformations.
 It was shown by Willmore \cite{Will} that, for surfaces in $ \R ^3$,  the Willmore energy satisfies the inequality $ \mc W ( f) \geq 4 \pi$ with equality if and only if the surface is a  round sphere. He conjectured also that  every  torus satisfies the inequality $\mc W(f) \geq 2 \pi^2$ with equality if and only if the surface is a M\"obius transform of the Clifford torus, what we will call \emph{conformal Clifford torus}. This conjecture was recently proved by Marques and Neves \cite{MN}. 

A natural question is the existence and energy values of non-minimizing critical points of the Willmore energy, the so called Willmore surfaces. Note that the Willmore energy is conformally invariant, therefore  inversions of complete non-compact minimal surfaces with appropriate growth at infinity are Willmore surfaces. It was proved by Bryant \cite{Bry} that smooth Willmore surfaces in $\R^3$ that are topologically spheres are all inversions of complete non-compact minimal spheres with embedded planar ends. Using the Weierstrass-Enneper representation, in principle, all these minimal surfaces may be classified. In particular, the Willmore energies of these surfaces are quantized, $ \mc W(f) = 4 \pi k, k \geq 1$ and the first non-trivial value is $\mc W(f) = 16 \pi$. Hence the gap to the next critical value of the Willmore energy among spheres is $ 12\pi$. The values $ k=2,3$ were recently investigated by Lamm and Nguyen \cite{LN} and it was shown they correspond to inversions of catenoids, Enneper's minimal surface or trinoids. These surfaces are not smooth but have point singularities or branch points.         

While, as just described, the family of Willmore spheres is quite well understood, much less is known for Willmore tori. By the work of Pinkall \cite{Pin}, it is known the existence of infinitely many Willmore tori  but it is an open problem whether or not the Willmore energies attained by critical tori  are isolated.
In this paper, we will show that there exists a gap from the minimizing conformal Clifford torus to the next critical point of the Willmore energy, namely we prove the following result,

\begin{theorem}[Gap Theorem for Willmore Tori]\label{Thm:Gap}
There exists $\epsilon_0>0$ such that, if $T \subset \R^3$ is a smoothly immersed Willmore torus (i.e. a critical point for the Willmore functional) with
\begin{equation}\label{gap}
\mc W(T)\leq 2 \pi^2 +\epsilon_0,
\end{equation}
then $T$ is the image,  under a M\"obius transformation of $\R^3$, of the standard Clifford torus $T_{Cl}$;  in particular $\mc W(T)=2 \pi^2$.
\end{theorem}

\begin{remark}
The analogous gap result for spheres in codimension one follows by the aforementioned work of Bryant \cite{Bry}, in codimension two was obtained by Montiel \cite{Mon}, and in arbitrary codimension was proved by  Kuwert and Sch\"atzle \cite{KSJDG} and subsequently by Bernard and Rivi\`ere \cite{BR}. 
\end{remark}
By the work of Chill-Fasangova-Sch\"atzle  \cite{CFS} (which is a continuation of  previous works on the Willmore flow by Kuwert and Sch\"atzle  \cite{KSJDG}, \cite{KS}) the gap Theorem \ref{Thm:Gap} has the following application to the Willmore  flow, i.e. the $L^2$-gradient flow of the Willmore functional:

\begin{corollary}\label{cor:WillFLow}
Let $ f _C:T \mapsto \R^3$ be a conformal Clifford torus. i.e.  a regular Moebius transformation of the Clifford torus $T_{Cl}$. 

There exists  $ \e_0>0$ such that if $ \| f_0 - f_C\| _{ W^{2,2}\cap C ^1 } \leq \e_0$ then, after reparametrisation by appropriate diffeomorphisms $ \Psi_t: T \rightarrow T$, 
the Willmore flow $(f_t)_t$ with initial data $ f_0$ exists globally and converges smoothly to  a conformal  Clifford torus $\tilde{f}_{C}$, that is 
\begin{align*}
f_t \circ \Psi_t \rightarrow \tilde{f}_C, \quad \text{ as $ t \rightarrow \infty $}.
\end{align*}     
\end{corollary}

We finish the paper by showing the following gap theorem from the Clifford torus to  higher genus (even non Willmore) surfaces. The proof is a combination of the large genus limit of the Willmore energy proved by Kuwert-Li-Sch\"atzle \cite{KLS} and the proof of the Willmore conjecture by Marques and Neves \cite{MN}.
 
\begin{theorem}\label{thm:HG}
There exists $ \e_0 > 0 $ such that if $\Sigma\subset \R ^ 3$ is a smooth immersed surface of genus at least two, then 
\begin{align*}
\mc W(\Sigma)\geq  2\pi ^ 2 + \e_0. 
\end{align*} 
\end{theorem} 

The  paper is organized as follows: in Section \ref{Sec:Prel} we gather together the necessary material and definitions that we will require in the paper.  In Section \ref{Sec:Gap} we prove  the gap Theorem \ref{Thm:Gap} for Willmore tori and, in Section \ref{Sec:Flow}, we use it to show the convergence of the Willmore flow of tori that are sufficiently close to the conformal Clifford torus, namely  we prove Corollary \ref{cor:WillFLow}. Finally,  Section \ref{sec:HG} is devoted to the proof of  the energy gap for higher genus surfaces, namely Theorem \ref{thm:HG}.

\section{Preliminaries} \label{Sec:Prel}
Throughout  the paper $T_{Cl}$ will denote the standard Clifford torus embedded in $\R^3$ (i.e. the torus obtained by revolution around the $z$-axis of a unit circle in the $xz$-plane with center at $(\sqrt{2},0,0)$) and $\CM$ will denote the M\"obius group of $\R^3$.

Recall from the classical paper of Weiner  (see in particular Lemma 3.3, Proposition 3.1 and Corollary 1 in \cite{Wei}), that  the second differential of the Willmore functional $\mc W''$ on the Clifford torus defines a positive semidefinite bounded symmetric bilinear form on $H^2(T_{Cl})$, the Sobolev space of $L^2$ integrable functions with first and second weak derivatives in $L^2$; in formulas, the second differential of the Willmore functional pulled back on $S^3$ via stereographic projection on the standard Clifford torus $\frac{1}{\sqrt{2}}(S^1\times S^1)\subset S^3\subset \R^4$ is given by
$$\mc W''(u,v):=\int (\Delta+|A|^2)[u]\cdot (\Delta+|A|^2+2)[v] \, d\mu_g,$$
where $\Delta$ is the Laplace-Beltrami operator on the Clifford torus, $|A|\equiv 2$ is the norm of the second fundamental form, and the integral is computed with respect to the surface measure $d\mu_g$.
\\ Moreover the kernel of the bilinear form consists of infinitesimal M\"obius transformations:
\begin{equation}\label{eq:Kernel}
K:=Ker(\mc W'')= \{\text{infinitesimal M\"obius transformations on } T_{Cl}\}\subset C^{\infty}(T_{Cl}),
\end{equation}
i.e. for every $w \in K$ there exists a M\"obius transformation $\Phi_w \in \CM$ such that for $t\in \R$ small enough
$$Exp_{T_{Cl}}(tw)=\Phi_{tw}(T_{Cl}),$$
where $Exp_{T_{Cl}}(tw)$ is the exponential in the normal direction with base surface $T_{Cl}$  and of magnitude $tw \in C^{\infty}(T_{Cl})$.
\\Again from \cite{Wei}, called $K^\perp\subset H^2(T_{Cl})$ the orthogonal space to $K$ in $H^2$, one also has that $\mc W''|_{K^{\perp}}$ is positive definite and defines a scalar product on $K^{\perp}$ equivalent to the restriction to $K^{\perp}$ of the $H^2(T_{Cl})$ scalar product:
\begin{equation}\label{eq:W''>0}
\mc W''(w,w)\geq \lambda \|w\|_{H^2(T_{Cl})}^2  \quad \forall\;  w \in K^{\perp},
\end{equation}
for some $\lambda>0$.

\section{Gap Theorem for Willmore Tori}\label{Sec:Gap}
The goal of this section is to prove Theorem \ref{Thm:Gap}.
\\

Recall from Section \ref{Sec:Prel} that $K\subset C^{\infty}(T_{Cl})\subset H^2(T_{Cl})$ denotes the kernel of the second differential $\mc W''$ of the Willmore functional on $T_{Cl}$ and it is made of infinitesimal M\"obius transformations; recall also that $K^{\perp}$ denotes the orthogonal complement of $K$ in $H^2(T_{Cl})$.
\\

Before proving the main theorem, let us prove a useful renormalization lemma  which roughly tells that if a surface $\Sigma$ is close to $T_{Cl}$ in $C^{k}$ topology, then up to a small conformal perturbation of $T_{Cl}$, $\Sigma$ can be written  as normal graph on $T_{Cl}$ via a function in $K^{\perp}\cap C^k$. 

\begin{lemma}\label{Lem:Normalise}
Let $T_{Cl}\subset \R^3$ be the standard Clifford torus and fix $2\leq k \in \N$. 

Then there exists a $\delta>0$ and continuous maps \footnote{notice that the map $u(\cdot)$ (respectively $v(\cdot)$ )  depends on another function $w$, so that $u(w)$ (resp. $v(w)$ ) denotes the value of the map $u$ (resp. $v$ ) evaluated at  $w$} 
\begin{eqnarray}
u&:&B^{C^k}_{\delta}(0) \to K \text{ endowed with } C^k \text {topology}, \quad \text{with } u(0)=0 \nonumber \\
v&:&B^{C^k}_{\delta}(0) \to C^k(T_{Cl})\cap K^{\perp},\quad \text{with } v(0)=0 \nonumber
\end{eqnarray}
where $B^{C^k}_{\delta}(0)\subset C^k(T_{Cl})$ is the ball of center $0$ and radius $\delta$,  such that the following holds:

For every $\Sigma\subset \R^3$ smooth closed (compact without boundary) embedded surface, $\delta-$close to $T_{Cl}$ in $C^{k}$ norm
\begin{equation}\label{Sigmaw}
\Sigma=Exp_{T_{Cl}}(w), \quad \text{for some } w \in C^k(T_{Cl}) \quad \text{with} \quad  \|w\|_{C^k(\Sigma)}< \delta,
\end{equation}
one has
\begin{equation}\label{eq:wuv}
\Sigma=Exp_{T_{Cl}}(w)=Exp_{[Exp_{T_{Cl}}(u(w))]}(v(w)).
\end{equation}
Notice that in the last term we are taking the exponential based on the surface $Exp_{T_{Cl}}(u(w))$ which is the image of $T_{Cl}$ under a small M\"obius transformation $\Phi_{u(w)}\in \CM$ since by construction $u(w)\in K$ is an infinitesimal M\"obius transformation (see Section \ref{Sec:Prel} for more details):
$$Exp_{T_{Cl}}(u(w))=\Phi_{u(w)}(T_{Cl}).$$
\end{lemma}

Notice that, for $\|w\|_{C^k}< \delta$ , we can identify (and we will always do it) the function spaces $C^k(\Sigma)$ and $C^k(T_{Cl})$ via the $C^k$ diffeomorphism from $T_{Cl}$ to $\Sigma$ induced by $Exp$.
\\

{\bf Proof of Lemma \ref{Lem:Normalise}}
\\Let $\Sigma=Exp_{T_{Cl}}(w)$ be as in \eqref{Sigmaw}. For $\delta>0$ small enough, it is clear that for $w,u \in C^k(T_{Cl})$ with $\|w\|_{C^k(T_{Cl})},\|u\|_{C^k(T_{Cl})} < \delta $ there exists a unique $v \in C^k(T_{Cl})$ such that
\begin{equation}\label{eq:wuv'}
\Sigma=Exp_{T_{Cl}}(w)=Exp_{[Exp_{T_{Cl}}(u)]}(v).
\end{equation}
It is also clear that, called $B^{C^k}_{\delta}(0)$ the ball of radius $\delta$ and center $0$ in $C^k(T_{Cl})$, the map
\begin{eqnarray}
\tilde{F}&:&B^{C^k}_{\delta}(0)\times (K\cap  B^{C^k}_{\delta}(0)) \to C^k(T_{Cl}) \nonumber \\
(w,u)&\mapsto & \tilde{F}(w,u):=v \text{ such that  \eqref{eq:wuv'} is satisfied} \nonumber 
\end{eqnarray}
is of class $C^1$. 
Denoted with $P_K:H^2(T_{Cl})\to K$ the orthogonal projection to the closed linear subspace $K=Ker(\mc W'')\subset H^2(T_{Cl})$, we define also
\begin{eqnarray}
F&:&B^{C^k}_{\delta}(0)\times (K\cap B^{C^k}_{\delta}(0)) \to K \nonumber \\
(w,u)&\mapsto & F(w,u):=P_K\circ \tilde{F}(w,u) \nonumber .
\end{eqnarray}
Since $\tilde{F}$ is $C^1$ then also $F$ is of class $C^1$.

Called $F_u(0,0)$ the $u$-partial derivative of $F$ computed in $(0,0)$, we have that
\begin{eqnarray}
F(0,0)&=&0 \nonumber \\
F_u(0,0)&=& -Id_{K}. \nonumber
\end{eqnarray}
By the Implicit Function Theorem  (see for instance Theorem 2.3 in \cite{AP}), taking $\delta>0$ maybe smaller, we conclude that there exists a $C^1$-function 
$$u(\cdot):B^{C^k}_{\delta}(0)\to K \text{ such that } u(0)=0 \text{ and } F(w,u(w))=0 \; \forall\, w \in B^{C^k}_{\delta}(0).$$ Therefore, the $C^1$ function $v(w):=\tilde{F}(w,u(w))$ satisfies 
$$P_K(v(w))=P_K\circ \tilde{F} (w,u(w))=F(w,u(w))=0\; \forall\, w \in B^{C^k}_{\delta}(0);$$
so $v\cdot):B^{C^k}_{\delta}(0)\to K^{\perp}$. 
\\The proof of the lemma is complete once we recall \eqref{eq:wuv'} and the definition of  $\tilde{F}$.   \hfill$\Box$
\\  

{\bf Proof of Theorem \ref{Thm:Gap}}
\\First of all, since by assumption $\mc W(T)<8\pi$, then $T$ is an embedded torus in $\R^3$ (this classical fact is proved in \cite{LY} or, by a monotonicity formula, in \cite{SiL}).

Let us prove the theorem by contradiction and assume that there exists a sequence $\{T_n\}_{n \in \N}$ of embedded Willmore tori in $\R^3$ such that 
\begin{equation}\label{eq:Wfk}
\mc W(T_n)\downarrow 2 \pi^2
\end{equation}
and $T_n$ is not the image of the Clifford torus up to any M\"obius transformation. 
\\

By the strong compactness of Willmore tori with energy strictly below $8 \pi$ (see Theorem 5.3 in \cite{KS} or Theorem I.8 in \cite{Riv1}), up to M\"obius transformations (which leave unchanged the Willmore functional, by the conformal invariance) and up to subsequences, we have that
$$T_n\to \tilde{T} \quad\text{smooth convergence of surfaces compactly contained  in } \R^3,$$
for some $\tilde{T}$ embedded torus in $\R^3$. By continuity of the Willmore functional under smooth convergence (actually  uniform $C^2$-convergence would be enough), by the assumption \eqref{eq:Wfk} we have that
\begin{equation}\label{eq:WTT}
\mc W(\tilde{T})=\lim_{n\uparrow \infty} \mc W(T_n)=2\pi^2. 
\end{equation}
From the recent proof of the Willmore conjecture by Marques and Neves (see Theorem A in \cite{MN}), $\tilde{T}$ must be the the image of the Clifford Torus $T_{Cl}$ under a M\"obius transformation of $\R^3$. Therefore, up to M\"obius transformations, for large $n$ the torus $T_n$ can be written as
\begin{equation}
T_n=Exp_{T_{Cl}}(w_n),\quad \|w_n\|_{C^k(T_{Cl})}\to  0 \text{ as } n \uparrow \infty \quad \forall k \in \N.
\end{equation}  
By Lemma \ref{Lem:Normalise}, taking $k=2$, we can write $T_n$ as
$$T_n=Exp_{Exp_{T_{Cl}}(u_n)} (v_n),$$
with $u_n\in K, v_n \in K^{\perp}$, $\|u_n\|_{C^2(T_{Cl})}+\|v_n\|_{C^2(T_{Cl})} \to 0$ as $n\uparrow \infty$. As remarked before, since $u_n \in K$ is an infinitesimal M\"obius transformation, there exist $\Phi_n \in \CM$ converging (in $C^2$-topology) to the identity as $n\uparrow \infty$ such that 
$$Exp_{T_{Cl}}(u_n)=\Phi_n(T_{Cl}).$$
Therefore, up to  M\"obius transformations, we can assume that 
$$T_n=Exp_{T_{Cl}}(\tilde{v}_n)$$
with $\|\tilde{v}_n\|_{C^2(T_{Cl})}\to 0$  and $\|P_K(\tilde{v}_n)\|_{C^2(T_{Cl})}=o(\|\tilde{v}_n\|_{C^2(T_{Cl})})$  as  $n\uparrow \infty$ (notice that the last remainder estimate comes combining the following facts: $\tilde{v}_n$ is obtained via pull back of $v_n$ with a M\"obius transformation, the space $K$ is made of infinitesimal M\"obius transformations, $P_K(v_n)=0$ by construction).
Since by the work of Weiner \cite{Wei} recalled in Section \ref{Sec:Prel} the  Willmore functional is twice differentiable on the Clifford torus $T_{Cl}$ with respect to $C^2$-normal variations and its first differential is given by the  fourth order nonlinear elliptic operator  \footnote{$A^\circ$ is the tracefree second fundamental form. i.e. $A^\circ=A-H g$}
$$\mc W'_{T_{Cl}}= \Delta H+ |A^\circ|^2 H \quad, $$
whose linearization (i.e. the second differential ${\mc W}''_{T_{Cl}}$) has the bilaplacian $\Delta^2$ as top order part, by a first oder expansion of $\mc W'$  we get that
$$\mc W'_{Exp_{T_{Cl}}(\varphi) }=\mc W'_{T_{Cl}}+\mc W''_{T_{Cl}}[\varphi]+o(\|\varphi\|_{C^2(T_{Cl})}) \quad \text{in } H^{-2}(T_{Cl}), $$
for every  $\varphi\in C^2(T_{Cl})$ with $\|\varphi\|_{C^2(T_{Cl})}$ small enough.
Taking now $\varphi=\tilde{v}_n$, we obtain
\begin{equation}\label{eq:W'W''}
\mc W'_{T_n}=\mc W'_{T_{Cl}}+\mc W''_{T_{Cl}}[\tilde{v}_n]+o(\|\tilde{v}_n\|_{C^2(T_{Cl})}) \quad \text{in } H^{-2}(T_{Cl})  ; 
\end{equation}
since the Clifford torus $T_{Cl}$ and the tori $T_n$, by assumption, are critical for $\mc W$, then the first two terms are identically null; therefore, called $\tilde{v}_n^{\perp}=\tilde{v}_n-P_K(\tilde{v}_n)\in K^\perp\cap C^2$ and recalled that $\|P_K(\tilde{v}_n)\|_{C^2(T_{Cl})}=o(\|\tilde{v}_n\|_{C^2(T_{Cl})})$, equation \eqref{eq:W'W''} yields
\begin{equation}\label{eq:quasi}
\|\mc W''_{T_{Cl}}[\tilde{v}_n^{\perp}]\|_{H^{-2}(T_{Cl})}=o(\|\tilde{v}^\perp_n\|_{C^2(T_{Cl})}). 
\end{equation}
Using \eqref{eq:W''>0}, since $\tilde{v}_n^{\perp} \in Ker(\mc W'')^{\perp}$, for some $\lambda>0$ we have
\begin{equation}\label{eq:quasi'}
\lambda \|\tilde{v}_n^{\perp}\|^2_{H^2(T_{Cl})}\leq \;\left<\mc W''_{T_{Cl}}[\tilde{v}_n^{\perp}],\tilde{v}_n^{\perp}  \right>_{H^{-2},H^2}= o(\|\tilde{v}^\perp_n\|^2_{C^2(T_{Cl})})
\end{equation}
where in the last equality we used \eqref{eq:quasi}.
\\Now recall that by construction $T_n=Exp_{T_{Cl}}(\tilde{v}_n)$ is a smooth Willmore torus with $\|\tilde{v}_n\|_{C^2}\to 0$, so by $\e$-regularity (see Theorem 2.10  in \cite{KSJDG} and notice that  for $n$ large enough, since $\|\tilde{v}_n\|_{C^1}$ is small then there exists $0<\delta<1$ independent of $n$ such that $(1-\delta)|A_{ij}|\leq |\partial^2_{ij} \tilde{v}_n| \leq (1+\delta) |A_{ij}|$ where $A_{ij}$ is the second fundamental form of $T_n$ and $\partial^2_{ij} \tilde{v}_n$ are the second derivatives of $\tilde{v}_n$; see also Theorem I.5 in \cite{Riv1}) we infer that there exists $C>0$ such that for $n$ large enough
$$
\|\tilde{v}_n\|_{C^2(T_{Cl})} \leq C \|\tilde{v}_n\|_{H^2(T_{Cl})}; 
$$
which, together with $\|P_K(\tilde{v}_n)\|_{C^2}=o(\|\tilde{v}_n\|_{C^2})$, gives that 
\begin{equation}\label{quasi''}
\|\tilde{v}_n^{\perp}\|_{C^2(T_{Cl})}\leq C \|\tilde{v}_n^{\perp}\|_{H^2(T_{Cl})}.
\end{equation}
Combining \eqref{eq:quasi'} and \eqref{quasi''} gives the contradiction
 $$\frac{\lambda}{C}\leq \frac {o(\|\tilde{v}_n^{\perp}\|_{C^2(T_{Cl})}^2)} {\|\tilde{v}_n^{\perp}\|_{C^2(T_{Cl})}^2} \to 0 \quad \text{as }\; n\uparrow \infty.$$
 \hfill$\Box$

\section{Willmore Flow of Tori} \label{Sec:Flow}
The Willmore flow is by definition the flow in the direction of the negative $L^2$ gradient of the Willmore functional.

For geometric flows, the short-time existence theory for smooth initial data is standard. The key fact is that although the equation is not strictly parabolic, the zeroes of the symbol of the differential operator are due only to the diffeomorphism invariance of the equation. By breaking this diffeomorphism invariance, the existence theory then reduces to standard parabolic theory. The long time existence and behavior of the flow is a much more delicate issue; indeed the flow can develop singularities in finite or infinite time or can converge to a stationary point of the functional. Here we are interested in the last possibility.

The Willmore flow has been the object of investigation of several papers by Kuwert and Sch\"atzle (see for instance  \cite{KSJDG} and \cite{KS}); in particular we will make use of the following result of Chill, Fasangova and   Sch\"atzle (see Theorem 1.2 in  \cite{CFS}).

\begin{theorem}[Chill-Fasangova-Sch\"atzle]\label{thm_CFS}
Let $ \Sigma$ be a closed surface and let $ f _W: \Sigma \to \R^3$ be a Willmore immersion that locally minimizes the Willmore functional in $ C^k (k\geq 2)$, in the sense that there exists $ \delta > 0$ such that
\begin{align*}
\mc W(f) \geq \mc W (f_W) \quad \text{whenever } \| f -f _W\|_{ C^k} \leq \delta \quad .
\end{align*} 

Then there exists $\e > 0$ such that for any immersion $f_0 : \Sigma \to \R^3$ satisfying $|f_0-f_W|_{W^{2,2}\cap C^1} < \e$, the corresponding Willmore flow $(f_t)_t$ with initial data $f_0$ exists globally and
converges smoothly, after reparametrization by appropriate diffeomorphisms $\Psi_t:\Sigma \to \Sigma$, to a Willmore immersion $f_\infty$ which also minimizes locally the
Willmore functional in $C^k$, i.e.
$$f_t \circ\Psi_t \to f_\infty \text{ as } t \to \infty.$$
\end{theorem}

This theorem applies to any surface $\Sigma_g$ that minimizes globally the Willmore functional among surfaces of genus $g$. In particular, thanks to  the proof of the Willmore conjecture by Marques and Neves \cite{MN}, this theorem applies to the conformal Clifford torus. 
\\

\textbf{Proof of Corollary \ref{cor:WillFLow}}
By \cite{MN}, any conformal Clifford torus $f_C:T\to \R^3$ globally minimizes the Willmore energy among tori, hence we can apply Theorem \ref{thm_CFS}.  Therefore  there exists $\e > 0$ such that for any immersion $f_0 : T \to \R^3$ satisfying $|f_0-f_C|_{W^{2,2}\cap C^1} < \e$, the corresponding Willmore flow $(f_t)_t$ with initial data $f_0$ exists globally and
converges smoothly, after reparametrization by appropriate diffeomorphisms $\Psi_t:T \to T$, to a Willmore immersion $f_\infty:T\to \R$. 
\\But, since the Willmore flow does not increase the Willmore energy, we have 
$$\mc W (f_\infty)\leq \mc W(f_0) \leq 2\pi^2+ \delta_\e, $$
for some $\delta_\e>0$ depending on $\e$ and such that $\lim_{\e\to 0} \delta_\e=0$.
\\For $\delta_\e \leq \e_0$ we can apply the gap Theorem \ref{Thm:Gap} and conclude that $f_\infty$ is a conformal Clifford torus. The claim follows.
\hfill$\Box$.

 \section{Gap Theorem for higher genus surfaces} \label{sec:HG}
In this section, we show that the Clifford torus is isolated  in Willmore energy with respect to higher genus (even non-Willmore) surfaces. 
\\

{\bf{Proof of Theorem \ref{thm:HG}}}.
\\For every $g\geq 2$,  Bauer and Kuwert \cite{BK} (inspired by the paper of Simon \cite{SiL} and by some ideas of Kusner \cite{Kus}; a  different proof of this result was given later by  Rivi\`ere \cite{Riv2}) proved that the infimum 
$$\beta_g:=\inf\{\mc W(\Sigma): \Sigma \subset \R^3 \text{ is a smoothly immersed genus $g$ surface}  \}$$
is attained by a smooth embedded surface $\Sigma_g \subset \R^3$ (notice that, even if one minimizes among \emph{immersed} surfaces, the minimizer is in fact \emph{embedded}. The reason being that by direct comparison with the stereographic projection of Lawson's genus $g$ minimal surfaces in $S^3$ one gets that $\beta_g <8\pi$.  Moreover by a result of Li-Yau \cite{LY} it is  known that every  closed surface $\Sigma$ smoothly immersed  in $\R^3$ with ${\mc W} (\Sigma)<8\pi$ is in fact embedded; so, once that  the existence of a smooth immersed minimizer for $\beta_g$  is proved, the embeddedness follows). Moreover, Kuwert, Li and Sch\"atzle \cite{KLS} proved that 
\begin{equation}\label{eq:limg}
\lim_{g \to \infty} \beta_g = 8 \pi.
\end{equation}
From the proof of the Willmore conjecture by Marques and Neves \cite{MN} we also know that 
\begin{equation}\label{eq:sigmag}
\mc W(\Sigma_g) > 2\pi^2, \quad \forall g>1.
\end{equation}
The claim follows then by the combination of  \eqref{eq:limg} and   \eqref{eq:sigmag}.

\hfill$\Box$

\end{document}